\documentclass[11pt]{article}
\usepackage{fullpage}
\usepackage{latexsym}
\usepackage{amsfonts}
\usepackage{amssymb}
\usepackage{amsmath}
\usepackage{color}
\usepackage{graphicx, graphics}
\usepackage{hyperref}
\usepackage{pstricks, pst-plot, epsfig}
\usepackage{xargs}
\usepackage{subcaption}
\usepackage{arydshln}
\newcommand{\ud}{\,\mathrm{d}}
\renewcommand\Re{\operatorname{Re}}
\renewcommand\Im{\operatorname{Im}}

\newcommand{\CC}{\mathbb{C}}

\begin{document}
%
%\title[Initial-to-Interface Maps]{Initial-to-Interface Maps for the Heat Equation on Composite Domains}
\title{Initial-to-Interface Maps for the Heat Equation on Composite Domains}
\date{\today}

\author{Natalie E Sheils\\School of Mathematics\\ University of Minnesota\\nesheils{@}umn.edu \and Bernard Deconinck\\Department of Applied Mathematics\\ University of Washington\\bernard{@}amath.washington.edu}
%\author[N.E. Sheils and B. Deconinck]{Natalie E. Sheils and Bernard Deconinck\thanks{nsheils{@}umn.edu,~deconic{@}uw.edu}}
%\affil{Department of Applied Mathematics\\
%University of Washington\\
%Seattle, WA 98195-2420
%}

\maketitle
\begin{center}\noindent \emph{This paper is dedicated to Mark Ablowitz on the occasion of his 70th birthday, in recognition of his many important contributions to nonlinear science.}
\end{center}
%%%%%%%%%%%%%%%%%%%%%%%%%%%%%%%%%%%%%%%%%%%%%%%%%%%

\begin{abstract}
A map from the initial conditions to the function and its first spatial derivative evaluated at the interface is constructed for the heat equation on finite and infinite domains with $n$ interfaces.  The existence of this map allows changing the problem at hand from an interface problem to a boundary value problem which allows for an alternative to the approach of finding a closed-form solution to the interface problem.  \end{abstract}
%%%%%%%%%%%%%%%%%%%%%%%%%%%%%%%%%%%%%%%%%%%%%%%%%%%

\section{Introduction}
Interface problems for partial differential equations (PDEs) are initial boundary value problems for which the solution of an equation in one domain prescribes boundary conditions for the equations in adjacent domains. In applications, conditions at the interface follow from conservations laws. Few interface problems allow for an explicit closed-form solution using classical solution methods. Using the Fokas method~\cite{FokasBook, FokasPelloni4} such solutions may be constructed.  This has been done in the case of the heat equation with $n$ interfaces in infinite, finite, and periodic domains as well as on graphs~\cite{Asvestas, DeconinckPelloniSheils, SheilsDeconinck_PeriodicHeat, SheilsSmith, Mantzavinos}.  The method has also been extended to dispersive problems~\cite{SheilsDeconinck_LS, SheilsDeconinck_LSp}, and higher order problems~\cite{DeconinckSheilsSmith}.  These works construct explicit solutions in terms of given initial and boundary conditions.  The value of the function at the interface is not known.  

In this paper we consider the heat equation with $n$ interfaces on domains of finite and infinite extent.  The problem of heat conduction in a composite wall is a classical problem in design and construction discussed in many excellent texts, see for instance~\cite{CarslawJaeger, HahnO}. It is usual to restrict to the case of walls whose constitutive parts are in perfect thermal contact and have physical properties that are constant throughout the material and that are considered to be of infinite extent in the directions parallel to the wall. Further, we assume that temperature and heat flux do not vary in these directions. In that case, the mathematical model for heat conduction in each wall layer is given by \cite[Chapter~10]{HahnO}:
\begin{subequations}
\begin{align}\label{heatgeneral}
&&u_t^{(j)}&=\alpha_j u_{xx}^{(j)}, & x_{j-1}<x<,x_j,\\
&&u^{(j)}(x,t=0)&=u_0^{(j)}(x), & x_{j-1}<x<,x_j,
\end{align}
\end{subequations}
here $u^{(j)}(x,t)$ denotes the temperature in the wall layer indexed by $(j)$, $\alpha_j>0$ is the heat-conduction coefficient of the $j$-th layer, $x=x_{j-1}$ is the left extent of the layer, and $x=x_j$ is its right extent. The sub-indices denote derivatives with respect to the one-dimensional spatial variable $x$ and the temporal variable $t$. The function $u^{(j)}_0(x)$ is the prescribed initial condition of the system. The continuity of the temperature $u^{(j)}$ and of its associated heat flux $\alpha_j u^{(j)}_x$ are imposed across the interface between layers. In what follows it is convenient to use the quantity $\sigma_j$, defined as the positive square root of $\alpha_j$: $\sigma_j=\sqrt{\alpha_j}$.

If the layer is either at the far left or far right of the wall, Dirichlet, Neumann, or Robin boundary conditions can be imposed on its far left or right boundary respectively, corresponding to prescribing ``outside'' temperature, heat flux, or a combination of these. A derivation of the interface boundary conditions is found in \cite[Chapter~1]{HahnO}. It should be noted that the set-up presented in \eqref{heatgeneral} also applies to the case of one-dimensional rods in thermal contact.  Even for the simple problem of two finite walls in thermal contact, the classical approach using separation of variables \cite{HahnO} can provide only an implicit answer.  Indeed, the solution obtained in~\cite{HahnO} depends on certain eigenvalues defined through a transcendental equation that can be solved only numerically. In contrast, the Fokas Method produces an explicit solution formula involving only known quantities. 

The construction of a Dirichlet-to-Neumann map, that is, determining the boundary values that are not prescribed in terms of the initial and boundary conditions, is important in the study of PDEs and particularly inverse problems~\cite{Fokas8, SylvesterUhlmann}.  In what follows we construct a similar map between the initial values of the PDE and the function (and some number of spatial derivatives) evaluated at the interface.  This map allows for an alternative to the approach of finding solutions to interface problems as presented in earlier papers using the Fokas method by the authors and others.  This would be most useful in the case where one is interested only in the behavior of solutions at the interface.  The method presented here can be extended in a straightforward way to many other interface problems.  To our knowledge, no such maps currently exist.  

Given the initial conditions, one could find the value of the function and its derivatives at the interface(s) using these maps.  This changes the problem at hand from an interface problem to a collection of independent boundary value problems (BVPs).  At this point, the BVPs could be solved using any number of methods appropriate for the given problem.  Each BVP would be over-specified, however, by construction it is clear that the corresponding spectral functions are admissible~\cite{FokasPelloni2001}, \emph{i.e.} the data is mutually compatible.

\section{The heat equation on an infinite domain with $n$ interfaces}\label{sec:I2I_nheat_i}
Consider
\begin{equation}\label{n_heat}
u_t=\sigma(x) u_{xx},
\end{equation}
together with the initial condition $u_0(x)=u(x,0)$ and the asymptotic conditions $\lim_{|x|\to\infty} u(x,t)=0$, where $-\infty<x<\infty$, $0<t<T$, and $$\sigma(x)=\left\{\begin{array}{lll} \sigma_1^2,&&x<x_1,\\ \sigma_2^2,&&x_1<x<x_2,\\\vdots \\
\sigma_n^2,&&x_{n-1}<x<x_n,\\
\sigma_{n+1}^2,&&x>x_n. \end{array}\right. $$  The restriction $\lim_{|x|\to\infty} u(x,t)=0$ can easily be made more general as in~\cite{DeconinckPelloniSheils}.

We can rewrite~\eqref{n_heat} as the set of equations
\begin{align}\label{I2I_nheat}
u^{(j)}_t=&\sigma_j^2 u^{(j)}_{xx}, &x_{j-1}<x<x_j, &~~0<t<T,
\end{align}
for $1\leq j\leq n+1$ where $x_0=-\infty$ and $x_{n+1}=\infty$.  We impose the continuity interface conditions~\cite{DeconinckPelloniSheils, Kevorkian}
\begin{align*}
u^{(j)}(x_j,t)=&u^{(j+1)}(x_j,t), &t>0,\\
\sigma_j^2u^{(j)}_x(x_j,t)=&\sigma_{j+1}^2u^{(j+1)}_x(x_j,t), &t>0,
\end{align*}
for $1\leq j\leq n$.  These interface conditions follow from conservation laws and are fully derived~\cite[Chapter~1]{SheilsThesis}. Since $u^{(j)}(x,t)$ is defined on the open interval $x_{j-1}<x<x_j$, when we write $u^{(j)}(x_j,t)$ we mean $\lim_{x\to x_j^-} u^{(j)}(x,t)$.  Similarly, we denote $\lim_{x\to x_j^+}u^{(j+1)}(x,t)$ by $u^{(j+1)}(x_j,t)$.  Without loss of generality we shift the problem so $x_1=0$.   Using the usual steps of the Fokas method~\cite{FokasBook, FokasPelloni4, DeconinckTrogdonVasan} we have the local relations

\begin{align}\label{I2I_heatn_local}
(e^{-ikx+\omega_j(k) t}u^{(j)}(x,t))_t=&(\sigma_j^2e^{-ikx+\omega_j(k) t}(u^{(j)}_x(x,t)+ik u^{(j)}(x,t)))_x,
\end{align}
where $\omega_j(k)=(\sigma_j k)^2$.  These relations are a one-parameter family obtained by rewriting~\eqref{I2I_nheat}.

\begin{figure}[htbp]
\begin{center}
\def\svgwidth{5in}
   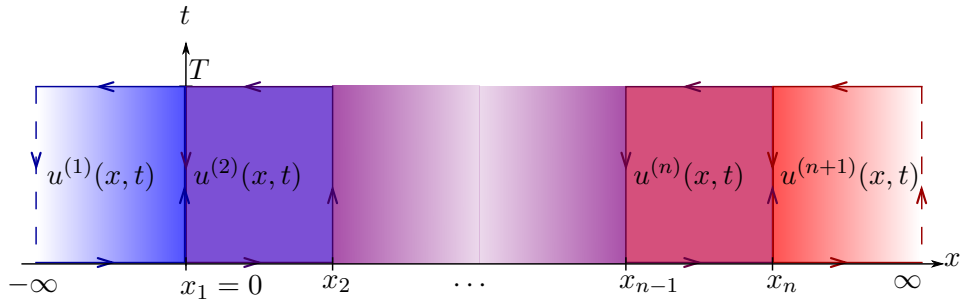 % requires the graphicx package
   \caption{Domains for the application of Green's Theorem in the case of an infinite domain with $n$ interfaces.   \label{fig:GR_domain_ni}}
  \end{center}
\end{figure}

\begin{figure}[htbp]
\begin{center}
\begin{subfigure}[b]{.45\textwidth}
\centering
\def\svgwidth{\textwidth}
   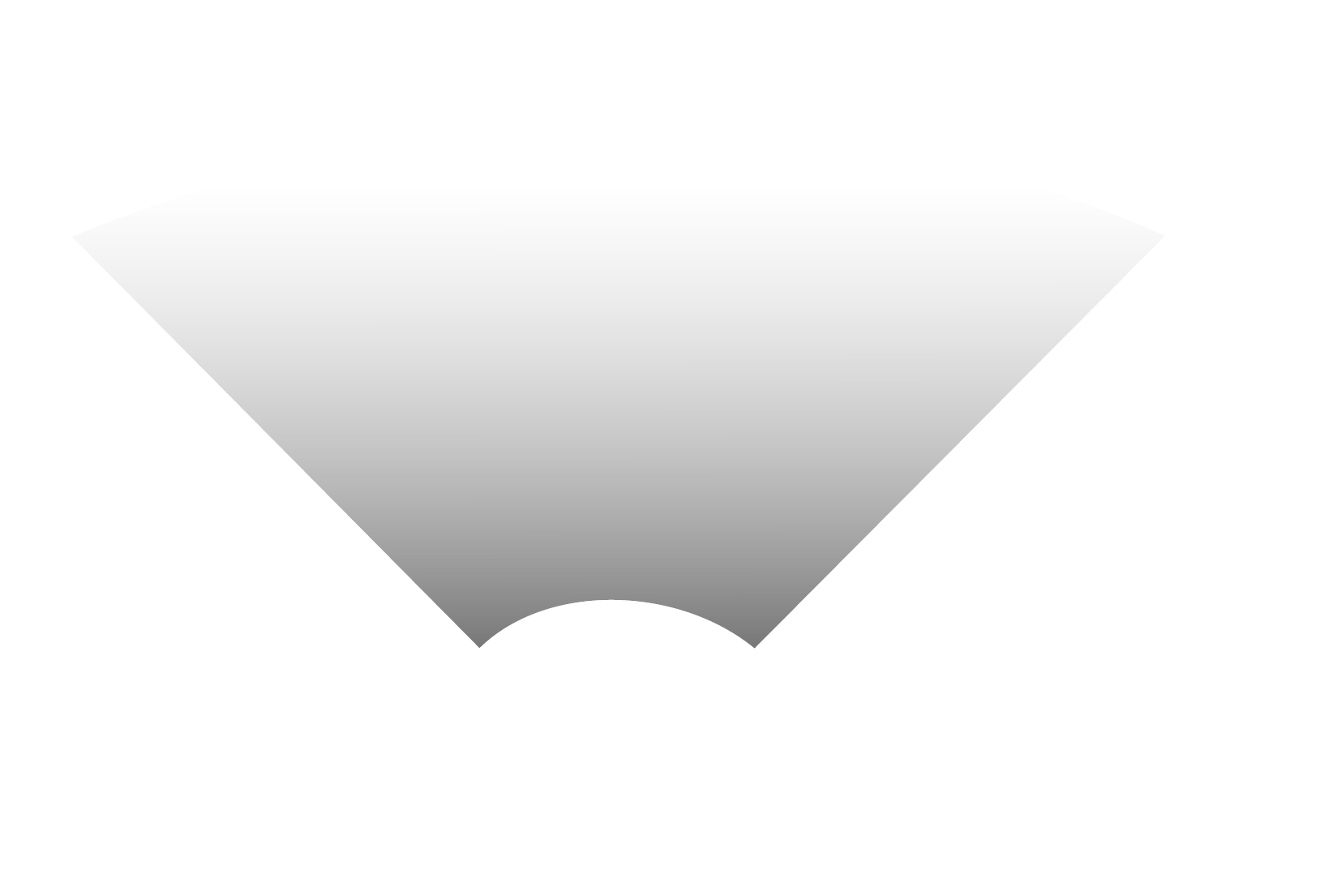 % requires the graphicx package
   \caption{ \label{fig:heat_DRpm}}
   \end{subfigure}
   \hfill
   \begin{subfigure}[b]{.45\textwidth}
\centering
\def\svgwidth{\textwidth}
   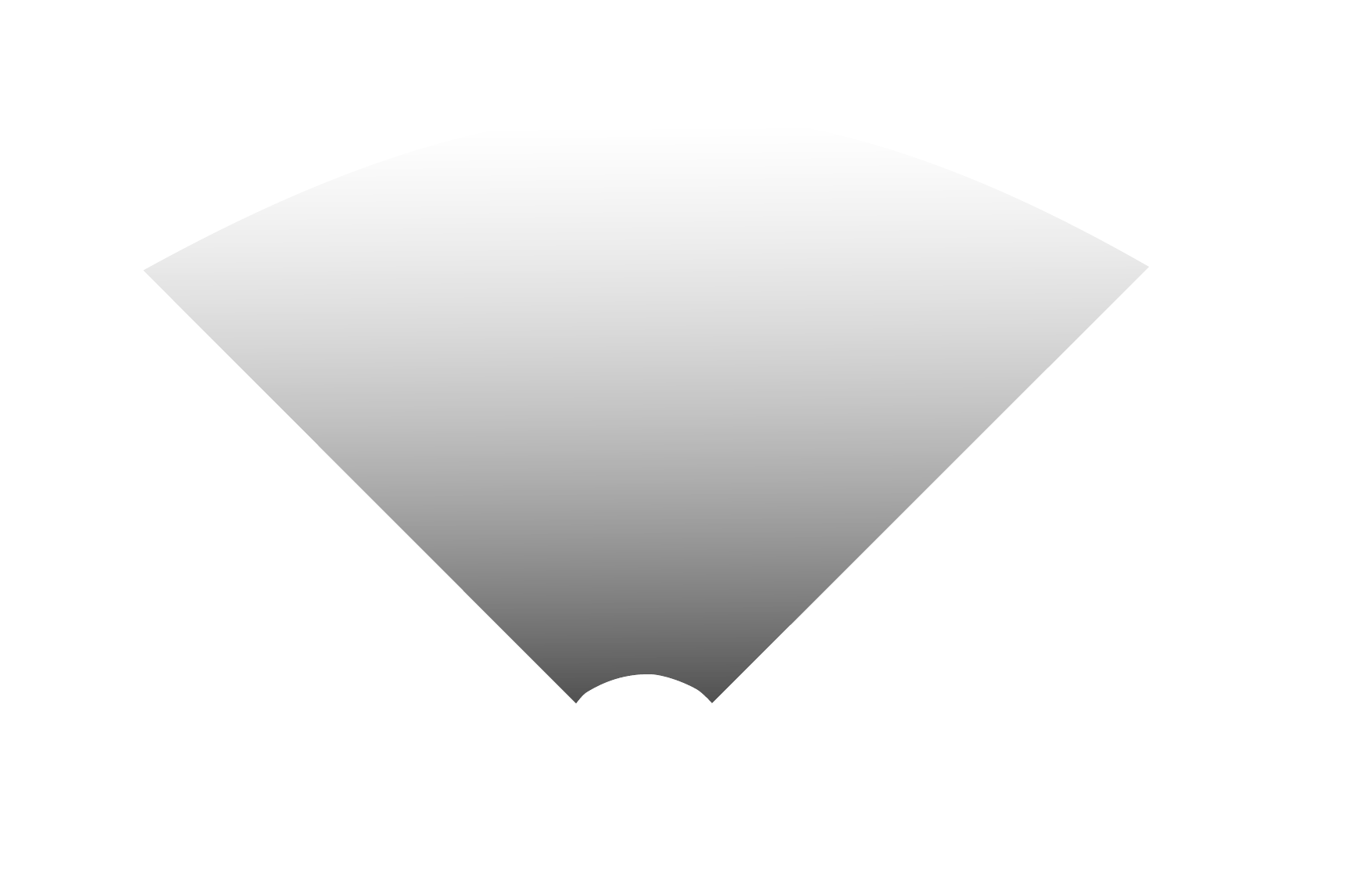 
   \caption{  \label{fig:heat_Dpm_close}}
\end{subfigure}
   \end{center}
   \caption{(a) The domain $D_R^+$ for the heat equation.  (b) The contour $\mathcal{L}^+$ is shown as a red dashed line.  An application of Cauchy's Integral Theorem using this contour allows elimination of the contribution of terms involving the Fourier transform of the solution. }
\end{figure}

Integrating over the appropriate cells of the domain (see Figure~\ref{fig:GR_domain_ni}) and applying Green's Theorem we find the global relations
\begin{equation}\label{nGR}
\begin{split}
0=&\int_{x_{j-1}}^{x_j} e^{-ikx}u^{(j)}_0(x)\ud x-\int_{x_{j-1}}^{x_j} e^{-ikx+\omega_j(k)T}u^{(j)}(x,T)\ud x\\
&+\int_0^T\sigma_j^2 e^{-ikx_j+\omega_j(k)s} (u^{(j)}_x(x_j,s)+ i k u^{(j)}(x_j,s))\ud s\\
&-\int_0^T\sigma_j^2 e^{-ikx_{j-1}+\omega_j(k)s} (u^{(j)}_x(x_{j-1},s)+ i k u^{(j)}(x_{j-1},s))\ud s,
\end{split}
\end{equation}
for $1\leq j\leq n+1$.  Define $D=\{k\in\CC: \Re(\omega_j(k))>0\}$, $D_R=\{k\in D: |k|>R\}$, and $D_R^+=\{k\in D_R: \Im(k)>0\}$  as in Figure~\ref{fig:heat_DRpm} where $R>0$ is an arbitrary finite constant. Since $|x|$ can become arbitrarily large for $j=1$ and $j=n+1$, we require $k\in\CC^+$ when $j=1$ and $k\in \CC^-$ when $j=n+1$, in Equation~\eqref{nGR} in order to guarantee that the integrals are defined. For $2\leq j\leq n$,~\eqref{nGR} is valid for $k\in\CC$.  The dispersion relation $\omega_j(k)=(\sigma_jk)^2$ is invariant under the symmetry $k\to -k$.  We supplement the $n+1$ global relations above with their evaluation at $-k$, namely,

\begin{equation}\label{nGR_minus}
\begin{split}
0=&\int_{x_{j-1}}^{x_j} e^{ikx}u^{(j)}_0(x)\ud x-\int_{x_{j-1}}^{x_j} e^{ikx+\omega_j(k)T}u^{(j)}(x,T)\ud x\\
&+\int_0^T\sigma_j^2 e^{ikx_j+\omega_j(k)s} (u^{(j)}_x(x_j,s)- i k u^{(j)}(x_j,s))\ud s\\
&-\int_0^T\sigma_j^2 e^{ikx_{j-1}+\omega_j(k)s} (u^{(j)}_x(x_{j-1},s)- i k u^{(j)}(x_{j-1},s))\ud s,
\end{split}
\end{equation}
for $1\leq j\leq n+1$.  When $j=1$,~\eqref{nGR_minus} is valid for $k\in\CC^-$. Similarly, for $j=n+1$,~\eqref{nGR_minus} is valid for $k\in\CC^+$.  For $2\leq j\leq n$,~\eqref{nGR_minus} is valid for all $k\in\CC$.  Without loss of generality we choose to work with the equations valid in the upper half plane.  Define
\begin{align*}
g^{(j)}_{0}({\omega},t)=&\int_{0}^te^{\omega s}u^{(j)}(x_j,s)\ud s=\int_{0}^te^{\omega s}u^{(j+1)}(x_j,s)\ud s,\\
g^{(j)}_{1}({\omega},t)=&\int_{0}^te^{\omega s}u_x^{(j)}(x_j,s)\ud s=\frac{\sigma_{j+1}^2}{\sigma_j^2}\int_{0}^te^{\omega s}u_x^{(j+1)}(x_j,s)\ud s,\\
\hat{u}^{(j)}(k,t)=&\int_{x_{j-1}}^{x_j}e^{-ikx}u^{(j)}(x,t)\ud x,\\
\hat{u}^{(j)}_0(k)=&\int_{x_{j-1}}^{x_j}e^{-ikx}u^{(j)}_0(x)\ud x,\\
\end{align*}
for $1\leq j\leq n$.  Using the change of variables $k=\kappa/\sigma_j$ on the $j^\textrm{th}$ equation, the global relations valid in the upper-half plane are
\begin{subequations}\label{heatn_GRu}
\begin{equation}
e^{\kappa^2 T}\hat{u}^{(1)}\left(\frac{\kappa}{\sigma_1},T\right)-\hat{u}_0^{(1)}\left(\frac{\kappa}{\sigma_1}\right)=e^{-i\kappa \frac{x_1}{\sigma_1}}\left(\frac{i\kappa}{\sigma_1} g_0^{(1)}(\kappa^2,T)+g_1^{(1)}(\kappa^2,T)\right),
\end{equation}
\begin{equation}
\begin{split}
e^{\kappa^2 T}\hat{u}^{(j)}\left(\frac{\kappa}{\sigma_j},T\right)-\hat{u}_0^{(j)}\left(\frac{\kappa}{\sigma_j}\right)=&e^{\frac{-i\kappa x_j}{\sigma_j}}\left(\frac{i\kappa}{\sigma_j} g_0^{(j)}(\kappa^2,T)+g_1^{(j)}(\kappa^2,T)\right)\\
&-e^{\frac{-i\kappa x_{j-1}}{\sigma_j}}\left(\frac{i\kappa}{\sigma_j} g_0^{(j-1)}(\kappa^2,T)+\frac{\sigma_{j-1}^2}{\sigma_j^2}g_1^{(j-1)}(\kappa^2,T)\right),
\end{split}
\end{equation}
\begin{equation}
\begin{split}
e^{\kappa^2 T}\hat{u}^{(j)}\left(\frac{-\kappa}{\sigma_j},T\right)-\hat{u}_0^{(j)}\left(\frac{-\kappa}{\sigma_j}\right)=&e^{\frac{i\kappa x_j}{\sigma_j}}\left(\frac{-i\kappa}{\sigma_j} g_0^{(j)}(\kappa^2,T)+g_1^{(j)}(\kappa^2,T)\right)\\
&+e^{\frac{i\kappa x_{j-1}}{\sigma_j}}\left(\frac{i\kappa}{\sigma_j} g_0^{(j-1)}(\kappa^2,T)-\frac{\sigma_{j-1}^2}{\sigma_j^2}g_1^{(j-1)}(\kappa^2,T)\right),
\end{split}
\end{equation}
\begin{equation}
e^{\kappa^2 T}\hat{u}^{(n+1)}\left(\frac{-\kappa}{\sigma_{n+1}},T\right)-\hat{u}_0^{(n+1)}\left(\frac{-\kappa}{\sigma_{n+1}}\right)=e^{\frac{i\kappa x_{n}}{\sigma_{n+1}}}\left(\frac{i\kappa}{\sigma_{n+1}} g_0^{(n)}(\kappa^2,T)-\frac{\sigma_{n}^2}{\sigma_{n+1}^2}g_1^{(n)}(\kappa^2,T)\right),
\end{equation}
\end{subequations}
for $2\leq j\leq n$.  Equation~\eqref{heatn_GRu} can be written as a linear system for the interface values:

$$\mathcal{A}(\kappa)X(\kappa^2,T)=Y(\kappa)+\mathcal{Y}(\kappa,T),$$
where
%\begin{subequations}
\begin{equation*}
X(\kappa^2,T)=\left(g_0^{(1)},g_0^{(2)},\ldots,g_0^{(n)},g_1^{(1)},g_1^{(2)},\ldots,g_1^{(n)}\right)^\top,
\end{equation*}
\begin{equation*}%\label{Yvector}
Y(\kappa)=-\left(\hat{u}_0^{(1)}\left(\frac{\kappa}{\sigma_1}\right),\ldots,\hat{u}_0^{(n)}\left(\frac{\kappa}{\sigma_n}\right),\hat{u}_0^{(2)}\left(\frac{-\kappa}{\sigma_2}\right),\ldots,\hat{u}_0^{(n+1)}\left(\frac{-\kappa}{\sigma_{n+1}}\right)\right)^\top,
\end{equation*}
\begin{equation*}
\mathcal{Y}(\kappa,T)=e^{\kappa^2T}\left(\hat{u}^{(1)}\left(\frac{\kappa}{\sigma_1},T\right),\ldots,\hat{u}^{(n)}\left(\frac{\kappa}{\sigma_n},T\right),\hat{u}^{(2)}\left(\frac{-\kappa}{\sigma_2},T\right),\ldots,\hat{u}^{(n+1)}\left(\frac{-\kappa}{\sigma_{n+1}},T\right)\right)^\top,
\end{equation*}
and

\begin{equation*}%\label{Amatrix}
\begin{split}
&\mathcal{A}(\kappa)=\\
&\left(
\begin{array}{lll:lll}
\frac{i\kappa}{\sigma_1}e^{-i\frac{\kappa x_1}{\sigma_1}}&&&e^{-i\frac{\kappa x_1}{\sigma_1}}\\
\frac{-i\kappa}{\sigma_2}e^{-i\frac{\kappa x_1}{\sigma_2}}&\frac{i\kappa}{\sigma_2}e^{-i\frac{\kappa x_2}{\sigma_2}}&&\frac{-\sigma_1^2}{\sigma_2^2}e^{-i\frac{\kappa x_1}{\sigma_2}}&e^{-i\frac{\kappa x_2}{\sigma_2}}\\
\hspace{.7in}\ddots&\hspace{.7in}\ddots&&\hspace{.7in}\ddots&\hspace{.7in}\ddots\\
&\frac{-i\kappa}{\sigma_n}e^{-i\frac{\kappa x_{n-1}}{\sigma_n}}&\frac{i\kappa}{\sigma_n}e^{-i\frac{\kappa x_n}{\sigma_n}}&&\frac{-\sigma_{n-1}^2}{\sigma_n^2}e^{-i\frac{\kappa x_{n-1}}{\sigma_n}}&e^{-i\frac{\kappa x_n}{\sigma_n}}\\
\hdashline
\frac{i\kappa}{\sigma_2}e^{i\frac{\kappa}{\sigma_2}x_1}&\frac{-i\kappa}{\sigma_2}e^{i\frac{\kappa x_2}{\sigma_2}}&&\frac{-\sigma_1^2}{\sigma_2^2}e^{i\frac{\kappa x_1}{\sigma_2}}&e^{i\frac{\kappa x_2}{\sigma_2}}\\
\hspace{.7in}\ddots&\hspace{.7in}\ddots&&\hspace{.7in}\ddots&\hspace{.7in}\ddots\\
&\frac{i\kappa}{\sigma_{n-1}}e^{i\frac{\kappa x_{n-1}}{\sigma_n}}&\frac{-i\kappa}{\sigma_n}e^{i\frac{\kappa x_n}{\sigma_n}}&&\frac{-\sigma_{n-1}^2}{\sigma_n^2}e^{i\frac{\kappa x_{n-1}}{\sigma_n}}&e^{i\frac{\kappa x_n}{\sigma_n}}\\
&&\frac{i\kappa}{\sigma_{n+1}}e^{i\frac{\kappa x_n}{\sigma_{n}}}&&&\frac{-\sigma_n^2}{\sigma_{n+1}^2}e^{i\frac{\kappa x_n}{\sigma_{n+1}}}
\end{array}
\right).
\end{split}
\end{equation*}
The matrix $\mathcal{A}(\kappa)$ consists of four $n\times n$ blocks as indicated by the dashed lines.  The two blocks in the upper half of $\mathcal{A}(\kappa)$ are zero except for entries on the main and $-1$ diagonals.  The lower two blocks of $\mathcal{A}(\kappa)$ only have nonzero entries on the main and $+1$ diagonals.  The matrix $\mathcal{A}(\kappa)$ is singular for isolated values of $\kappa$.  Asymptotically, for large $|\kappa|$, the zeros of $\det(\mathcal{A}(\kappa))$ lie within a strip parallel to the real line~\cite{Langer}.  Since asymptotically there are no zeros in $D_R^+$, a sufficiently large $R$ may be chosen such that $\mathcal{A}(\kappa)$ is nonsingular for every $\kappa\in D_R^+$ and $\det(\mathcal{A}(\kappa))\neq0$. 

%\medskip
%\textbf{Remark.}  We have been unable to construct physical examples where the zeros of $\det(\mathcal{A}(\kappa))$ are in $D^+$ and are different from 0.  However, if nonphysical values of the parameters are chosen (\emph{e.g.}, $\sigma_j$ imaginary), then  $\det(\mathcal{A}(\kappa))$ may zeros in $D^+$. 
%
%\textcolor{red}{What does one do then? If the zeros are finitely many, is it meaningful to deform the contours in first place so that these zeros are eventually not enclosed in $D^+$? }
%\medskip

Using Cramer's Rule to solve this system, we have
\begin{subequations}\label{I2I_g01solns}
\begin{align}
g_0^{(j)}(\kappa^2,T)=&\frac{\det(\mathcal{A}_j(\kappa,T))}{\det(\mathcal{A}(\kappa))},\\
g_1^{(j)}(\kappa^2,T)=&\frac{\det(\mathcal{A}_{j+n}(\kappa,T))}{\det(\mathcal{A}(\kappa))},
\end{align}
\end{subequations}
where $1\leq j\leq n$ and $\mathcal{A}_j(\kappa,T)$ is the matrix $\mathcal{A}(\kappa)$ with the $j^\textrm{th}$ column replaced by $Y+\mathcal{Y}$.  This does not give an effective initial-to-interface map because~\eqref{I2I_g01solns} depends on the solutions $\hat{u}^{(j)}(\cdot,T)$.  To eliminate this dependence we multiply~\eqref{I2I_g01solns} by $\kappa e^{-\kappa^2t}$ and integrate around $D_R^+$, as is typical in the construction of Dirichlet-to-Neumann maps~\cite{FokasBook}. Switching the order of integration we have
\begin{subequations}\label{gsolns}
\begin{align}
\int_0^T u^{(j)}(x_j,s) \int_{\partial D_R^+}\kappa e^{\kappa^2(s-t)}\ud \kappa\ud s=&\int_{\partial D_R^+} e^{-\kappa^2t} \frac{\kappa\det(\mathcal{A}_j(\kappa,T))}{\det(\mathcal{A}(\kappa))}\ud \kappa, \\
\int_0^T u^{(j)}_x(x_j,s) \int_{\partial D_R^+} \kappa e^{\kappa^2(s-t)}\ud \kappa\ud s=&\int_{\partial D_R^+} e^{-\kappa^2t}\frac{\kappa \det(\mathcal{A}_{j+n}(\kappa,T))}{\det(\mathcal{A}(\kappa))}\ud \kappa.
\end{align}
\end{subequations}
Using the change of variables $i\ell=\kappa^2$ and the classical Fourier transform formula for the delta function we have
\begin{subequations}\label{I2I_INTg01solns}
\begin{align}
u^{(j)}(x_j,t) =&\frac{1}{i\pi}\int_{\partial D_R^+} e^{-\kappa^2t} \frac{\kappa\det(\mathcal{A}_j(\kappa,T))}{\det(\mathcal{A}(\kappa))}\ud \kappa, \\
u^{(j)}_x(x_j,t)=&\frac{1}{i\pi}\int_{\partial D_R^+} e^{-\kappa^2t}\frac{\kappa \det(\mathcal{A}_{j+n}(\kappa,T))}{\det(\mathcal{A}(\kappa))}\ud \kappa.
\end{align}
\end{subequations}

To examine the right-hand-side of~\eqref{I2I_INTg01solns} we factor the matrix $\mathcal{A}(\kappa)$ as $\mathcal{A}^L(\kappa)\mathcal{A}^M(\kappa)$ where
\begin{equation*}
\mathcal{A}^L(\kappa)=\left(
\begin{array}{lll:lll}
e^{-i\frac{\kappa}{\sigma_1}x_1}&&&\\
&e^{-i\frac{\kappa}{\sigma_2}x_2}&&&\\
&\hspace{.7in}\ddots&&\\
&&e^{-i\frac{\kappa}{\sigma_n}x_n}&&&\\
\hdashline
&&&e^{i\frac{\kappa}{\sigma_2}x_1}&\\
&&&&e^{i\frac{\kappa}{\sigma_3}x_2}&\\
&&&&\hspace{.7in}\ddots\\
&&&&&e^{i\frac{\kappa}{\sigma_{n+1}}x_n}
\end{array}
\right)
\end{equation*}
is a diagonal matrix.  The elements of $\mathcal{A}^M(\kappa)$ are either $0$, $\mathcal{O}(\kappa)$, or decaying exponentially fast for $\kappa\in D_R^+$.  Hence, $$\det(\mathcal{A}^M(\kappa))=c(\kappa)=\mathcal{O}(\kappa^{2n}),$$ for large $\kappa$ in $D_R^+$.  Now, $\det(\mathcal{A}(\kappa))=c(\kappa) \det(\mathcal{A}^L(\kappa))$ as $\kappa\to\infty$ for $\kappa\in D_R^+$.  Similarly, factor $\mathcal{A}_j(\kappa,T)=\mathcal{A}^L(\kappa) \mathcal{A}_j^M(\kappa,T)\mathcal{A}_j^R(\kappa,T)$ where $\mathcal{A}_j^R(\kappa,T)$ is the $2n\times 2n$ identity matrix with the $(j,j)$ component replaced by $e^{\kappa^2 T}$.   Then $\det(\mathcal{A}_j(\kappa,T))= e^{\kappa^2 T} \det(\mathcal{A}^L(\kappa)) \det(\mathcal{A}_j^M(\kappa,T))$.  Thus, the integrand we are considering in~\eqref{I2I_INTg01solns} is

$$\int_{\partial D_R^+}e^{-\kappa^2 t} \frac{\kappa \det(\mathcal{A}_j(\kappa,T))}{\det(\mathcal{A})}\ud \kappa=\int_{\partial D_R^+}e^{\kappa^2 (T-t)} \frac{\kappa \det(\mathcal{A}_j^M(\kappa,T))}{c(\kappa)}\ud \kappa.$$

The elements of $\mathcal{A}_j^M(\kappa,T)$ are the same as those in $\mathcal{A}^M(\kappa)$ except in the $j^\textrm{th}$ column.  Expanding the determinant of $\mathcal{A}_j^M(\kappa,T)$ along the $j^\textrm{th}$ column we see that 
\begin{align}\label{badstuff}
\begin{split}
e^{\kappa^2 (T-t)}&\frac{\kappa\det(\mathcal{A}_j^M(\kappa,T))}{c(\kappa)}=\sum_{\ell=1}^n\left( c_\ell(\kappa) \left(e^{\frac{i\kappa x_\ell}{\sigma_\ell}+\kappa^2(T-t)}\hat{u}^{(\ell)}\left(\frac{\kappa}{\sigma_\ell},T\right)-e^{-\kappa^2 t +\frac{i\kappa x_\ell}{\sigma_\ell}}\hat{u}_0^{(\ell)}\left(\frac{\kappa}{\sigma_\ell}\right)\right)\right.\\
&+\left. c_{\ell+n}(\kappa) \left(e^{\frac{-i\kappa x_\ell}{\sigma_{\ell+1}}+\kappa^2(T-t)}\hat{u}^{(\ell+1)}\left(\frac{-\kappa}{\sigma_{\ell+1}},T\right)-e^{-\kappa^2 t -\frac{i\kappa x_\ell}{\sigma_{\ell+1}}}\hat{u}_0^{(\ell+1)}\left(\frac{-\kappa}{\sigma_{\ell+1}}\right)\right)\right),
\end{split}
\end{align}
where $c_{\ell}=\mathcal{O}(\kappa^{0})$ and $c_{\ell+n}=\mathcal{O}(\kappa)$ for $1\leq\ell\leq n$.  The terms involving $\hat{u}^{(\ell)}(\cdot,T)$, the solutions of our equation, are decaying exponentially for $\kappa\in D_R^+$.  Thus, by Jordan's Lemma~\cite{AblowitzFokas}, the integral of this term along a closed, bounded curve in $\CC^+$ vanishes. In particular we consider the closed curve $\mathcal{L}^+=\mathcal{L}_{D_R^+}\cup\mathcal{L}^+_C$ where $\mathcal{L}_{D_R^+}=\partial D_R^+ \cap \{k: |k|<C\}$ and $\mathcal{L}^+_C=\{k\in D_R^+: |k|=C\}$, see Figure~\ref{fig:heat_Dpm_close}. Since the integral along $\mathcal{L}_C^+$ vanishes for large $C$,~\eqref{badstuff} must vanish since the contour $\mathcal{L}_{D_R^+}$  becomes $\partial D_R^+$ as $C\to\infty$.  

Since the terms involving the elements of $\mathcal{Y}(\kappa,T)$ evaluate to zero in the solution expression we have the solution
\begin{subequations}\label{I2I_INTg01solns_final}
\begin{align}
u^{(j)}(x_j,t) =&\frac{1}{i\pi}\int_{\partial D_R^+} e^{-\kappa^2t} \frac{\kappa\det(A_j(\kappa))}{\det(\mathcal{A}(\kappa))}\ud \kappa, \\
u^{(j)}_x(x_j,t)=&\frac{1}{i\pi}\int_{\partial D_R^+} e^{-\kappa^2t}\frac{\kappa \det(A_{j+n}(\kappa))}{\det(\mathcal{A}(\kappa))}\ud \kappa,
\end{align}
\end{subequations}
for $1\leq j\leq n+1$, where $A_j(\kappa)$ is the matrix $\mathcal{A}(\kappa)$ with the $j^\textrm{th}$ column replaced by $Y(\kappa)$. Equation~\ref{I2I_INTg01solns_final} is an effective map between the values of the function at the interface and the given initial conditions.  

\medskip
\textbf{Remark.} Note that since the problem is linear, one could have assumed the initial condition was zero for $x$ outside the region $x_{\ell-1}<x<x_\ell$.  Then, the map would be in terms of just $u_0^{(\ell)}(\cdot)$.  Summing over $1\leq \ell\leq n+1$ would give the complete map for a general initial condition. 
\medskip

As an example of a specific initial-to-interface map we consider the equation~\eqref{n_heat} with $n=1$.  Using~\eqref{I2I_INTg01solns_final} we have
\begin{align*}
\sigma_1^2u^{(1)}_x(0,t)=&\frac{i \sigma_1\sigma_2}{\pi(\sigma_1+\sigma_2)}\int_{\partial D_R^+}\kappa e^{-\kappa^2t}\left(\sigma_1\hat{u}_0^{(1)}\left(\frac{\kappa}{\sigma_1}\right)-\sigma_2\hat{u}_0^{(2)}\left(\frac{-\kappa}{\sigma_2}\right)\right)\ud \kappa,\\
u^{(1)}(0,t)=&\frac{1}{\pi(\sigma_1+\sigma_2)}\int_{\partial D_R^+}  e^{-\kappa^2t}\left(\sigma_1^2\hat{u}_0^{(1)}\left(\frac{\kappa}{\sigma_1}\right)+\sigma_2^2\hat{u}_0^{(2)}\left(\frac{-\kappa}{\sigma_2}\right)\right)\ud \kappa.
\end{align*}
In this case we can deform $D_R^+$ back to the real line with no pole contributions.  %For general $n$ this is not the case (see the following Remark).  
Switching the order of integration and evaluating the $\kappa$ integral we have
\begin{subequations}\label{1heat}
\begin{align}
\sigma_1^2u^{(1)}_x(0,t)=&\frac{\sigma_1\sigma_2}{2t^{3/2}\sqrt{\pi}(\sigma_1+\sigma_2)}\left(\int_{-\infty}^0 ye^{\frac{-y^2}{4t\sigma_1^2}}u_0^{(1)}(y)\ud y+\int_{0}^\infty ye^{\frac{-y^2}{4t\sigma_2^2}}u_0^{(2)}(y)\ud y\right),\\
u^{(1)}(0,t)=&\frac{1}{\sqrt{\pi t}(\sigma_1+\sigma_2)}\left(\sigma_1^2 \int_{-\infty}^0 e^{\frac{-y^2}{4t\sigma_1^2}}u_0^{(1)}(y)\ud y+\sigma_2^2 \int_{0}^\infty e^{\frac{-y^2}{4t\sigma_2^2}}u_0^{(2)}(y)\ud y\right),
\end{align}
\end{subequations}
which is an explicit map from the initial data to the value of the temperature and its associated flux at the interface, $x=0$.  If one allows $\sigma_1=\sigma_2$ the problem is simply that of the heat equation on the whole line. Equation~\eqref{1heat} with $\sigma_1=\sigma_2$ is exactly the Green's Function solution of the whole line problem evaluated at $x=0$~\cite{Kevorkian}. 

%\medskip
%\textbf{Remark.} Already in the case $n=2$, $\det(\mathcal{A}(\kappa))$ has zeros whenever $$\kappa=\frac{\sigma_2}{2x_2}\left(2\pi m -i \log\left(\frac{-(\sigma_1+\sigma_2)(\sigma_2+\sigma_3)}{(\sigma_1-\sigma_2)(\sigma_2-\sigma_3)}\right)\right),$$ for any integer $m$.  Thus, deforming $\partial D_R^+$ in~\eqref{I2I_INTg01solns_final} would result in an integral along the real line as well as a sum of residues from the contributions at the poles.  For larger $n$, it may not be possible to find explicitly the location of the poles.
 
\section{The heat equation on a finite domain with $n$ interfaces}\label{sec:I2I_nheat_f}

Consider~\eqref{n_heat} on a finite domain,  $x_0\leq x\leq x_{n+1}$, with the boundary conditions
\begin{subequations}\label{I2I_bcs}
\begin{align}
\beta_1 u^{(1)} (x_0 ,t)+\beta_2 u^{(1)} _x(x_0 ,t)=&f_1  (t),&t>0,\\
\beta_3 u^{(n+1)} (x_{n+1} ,t)+\beta_4  u^{(n+1)} _x(x_{n+1} ,t)=&f_2  (t),&t>0.
\end{align}
\end{subequations}
As before, we rewrite~\eqref{n_heat} as the set of equations
\begin{align*}%\label{I2I_nheat_f}
u^{(j)}_t=&\sigma_j^2 u^{(j)}_{xx}, &x_{j-1}<x<x_j, &~~0<t<T,
\end{align*}
for $1\leq j\leq n+1$, subject to the continuity interface conditions
\begin{align*}
u^{(j)}(x_j,t)=&u^{(j+1)}(x_j,t), &t>0,\\
\sigma_j^2u^{(j)}_x(x_j,t)=&\sigma_{j+1}^2u^{(j+1)}_x(x_j,t), &t>0,
\end{align*}
for $1\leq j\leq n$.  Without loss of generality we shift the problem so that $x_0=0$. 

\begin{figure}[htbp]
\begin{center}
\def\svgwidth{5in}
   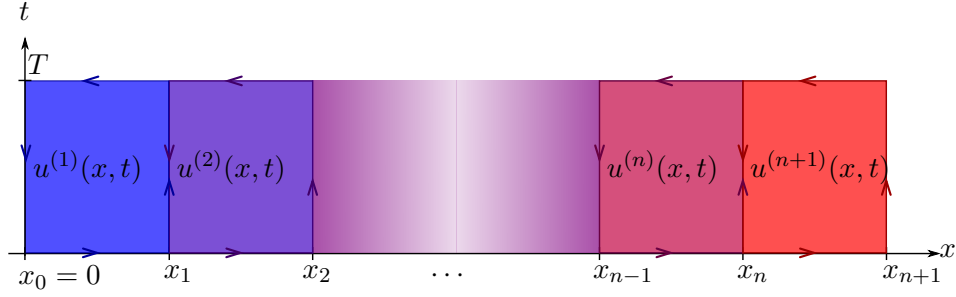 % requires the graphicx package
   \caption{Domains for the application of Green's Theorem in the case of a finite domain with $n$ interfaces.   \label{fig:GR_domain_nf}}
  \end{center}
\end{figure}
The following steps are very similar to those presented in the previous section.  In what follows we give a brief outline of the changes needed to solve on a finite domain.

Integrating the local relations~\eqref{I2I_heatn_local} around the appropriate domain (see Figure~\ref{fig:GR_domain_ni}) and applying Green's Theorem we find the global relations~\eqref{nGR} and their evaluation at $-k$~\eqref{nGR_minus}.  In contrast to Section~\ref{sec:I2I_nheat_i}, these $2n+2$ global relations are all valid for $k\in\CC$. 
%\begin{equation}\label{nGR_f}
%\begin{split}
%0=&\int_{x_{j-1}}^{x_j} e^{-ikx}u^{(j)}_0(x)\ud x-\int_{x_{j-1}}^{x_j} e^{-ikx+\omega_j(k)T}u^{(j)}(x,T)\ud x\\
%&+\int_0^T\sigma_j^2 e^{-ikx_j+\omega_j(k)s} (u^{(j)}_x(x_j,s)+ i k u^{(j)}(x_j,s))\ud s\\
%&-\int_0^T\sigma_j^2 e^{-ikx_{j-1}+\omega_j(k)s} (u^{(j)}_x(x_{j-1},s)+ i k u^{(j)}(x_{j-1},s))\ud s,
%\end{split}
%\end{equation}
%for $1\leq j\leq n+1$.  
%The dispersion relation $\omega_j(k)=(\sigma_jk)^2$ is invariant under the symmetry $k\to -k$.  We supplement the $n+1$ global relations above with their evaluation at $-k$, namely,
%\begin{equation}\label{nGR_minus_f}
%\begin{split}
%0=&\int_{x_{j-1}}^{x_j} e^{ikx}u^{(j)}_0(x)\ud x-\int_{x_{j-1}}^{x_j} e^{ikx+\omega_j(k)T}u^{(j)}(x,T)\ud x\\
%&+\int_0^T\sigma_j^2 e^{ikx_j+\omega_j(k)s} (u^{(j)}_x(x_j,s)- i k u^{(j)}(x_j,s))\ud s\\
%&-\int_0^T\sigma_j^2 e^{ikx_{j-1}+\omega_j(k)s} (u^{(j)}_x(x_{j-1},s)- i k u^{(j)}(x_{j-1},s))\ud s,
%\end{split}
%\end{equation}
%for $1\leq j\leq n+1$ which is again valid for $k\in\CC\setminus D$.  
In addition to the definitions in Section~\ref{sec:I2I_nheat_i} we define
\begin{align*}
g^{(0)}_{0}({\omega},t)=&\int_{0}^te^{\omega s}u^{(1)}(x_0,s)\ud s,\\
g^{(n+1)}_{0}({\omega},t)=&\int_{0}^te^{\omega s}u^{(n+1)}(x_{n+1},s)\ud s,\\
g^{(0)}_{1}({\omega},t)=&\int_{0}^te^{\omega s}u_x^{(1)}(x_0,s)\ud s,\\
g^{(n+1)}_{1}({\omega},t)=&\int_{0}^te^{\omega s}u_x^{(n+1)}(x_{n+1},s)\ud s,\\
\tilde{f}_m(\omega,t)=&\int_0^t e^{\omega s}f_m(s)\ud s,\\
\end{align*}
for $m=1,2$.  Using the change of variables $k=\kappa/\sigma_j$, the $n+1$ global relations are
\begin{subequations}
\begin{equation}
\begin{split}
e^{\kappa^2 t}\hat{u}^{(j)}\left(\frac{\kappa}{\sigma_j},T\right)&-\hat{u}_0^{(j)}\left(\frac{\kappa}{\sigma_j}\right)=e^{\frac{-i\kappa x_j}{\sigma_j}}\left(\frac{i\kappa}{\sigma_j} g_0^{(j)}(\kappa^2,T)+g_1^{(j)}(\kappa^2,T)\right)\\
&-e^{\frac{-i\kappa x_{j-1}}{\sigma_j}}\left(\frac{i\kappa}{\sigma_j} g_0^{(j-1)}(\kappa^2,T)+\frac{\sigma_{j-1}^2}{\sigma_j^2}g_1^{(j-1)}(\kappa^2,T)\right),
\end{split}
\end{equation}
\begin{equation}
\begin{split}
e^{\kappa^2 t}\hat{u}^{(j)}\left(\frac{-\kappa}{\sigma_j},T\right)&-\hat{u}_0^{(j)}\left(\frac{-\kappa}{\sigma_j}\right)=e^{\frac{i\kappa x_j}{\sigma_j}}\left(\frac{-i\kappa}{\sigma_j} g_0^{(j)}(\kappa^2,T)+g_1^{(j)}(\kappa^2,T)\right)\\
&+e^{\frac{i\kappa x_{j-1}}{\sigma_j}}\left(\frac{i\kappa}{\sigma_j} g_0^{(j-1)}(\kappa^2,T)-\frac{\sigma_{j-1}^2}{\sigma_j^2}g_1^{(j-1)}(\kappa^2,T)\right),
\end{split}
\end{equation}
\end{subequations}
for $1\leq j\leq n+1$ where we define $\sigma_0=\sigma_1$ for convenience.  These equations, together with the boundary values~\eqref{I2I_bcs}, can be written as a linear system for the interface values

$$\mathcal{A}^F X^F=Y^F+\mathcal{Y}^F,$$
where
\begin{subequations}
{\small
\begin{equation*}
X^F(\kappa^2,T)=\left(g_0^{(0)},g_0^{(1)},\ldots,g_0^{(n+1)},g_1^{(0)},g_1^{(1)},\ldots,g_1^{(n+1)}\right)^\top,
\end{equation*}
\begin{equation*}
Y^F(\kappa,T)=-\left(-\tilde{f}_1(\kappa^2, T),\hat{u}_0^{(1)}\left(\frac{\kappa}{\sigma_1}\right),\ldots,\hat{u}_0^{(n+1)}\left(\frac{\kappa}{\sigma_n}\right),\hat{u}_0^{(1)}\left(\frac{-\kappa}{\sigma_1}\right),\ldots,\hat{u}_0^{(n+1)}\left(\frac{-\kappa}{\sigma_{n+1}}\right),-\tilde{f}_2(\kappa^2, T)\right)^\top,
\end{equation*}
\begin{equation*}
\mathcal{Y}^F(\kappa,T)=e^{\kappa^2T}\left(0,\hat{u}^{(1)}\left(\frac{\kappa}{\sigma_1},T\right),\ldots,\hat{u}^{(n+1)}\left(\frac{\kappa}{\sigma_n},T\right),\hat{u}^{(1)}\left(\frac{-\kappa}{\sigma_1},T\right),\ldots,\hat{u}^{(n+1)}\left(\frac{-\kappa}{\sigma_{n+1}},T\right),0\right)^\top,
\end{equation*}
}
and
{\small
\begin{align*}\label{Amatrix_f}
&\mathcal{A}^F(\kappa)=\\
&\left(
\begin{array}{ccc:ccc}
\beta_1&&&\beta_2\\
\frac{-i\kappa}{\sigma_1}e^{-i\frac{\kappa x_0}{\sigma_1}}&\frac{i\kappa}{\sigma_1}e^{-i\frac{\kappa x_1}{\sigma_1}}&&-\frac{\sigma_0^2}{\sigma_1^2}e^{-i\frac{\kappa x_0}{\sigma_1}}&e^{-i\frac{\kappa x_1}{\sigma_1}}\\
\hspace{.7in}\ddots&\hspace{.7in}\ddots&&\hspace{.7in}\ddots&\hspace{.7in}\ddots\\
&\frac{-i\kappa}{\sigma_{n+1}}e^{-i\frac{\kappa x_{n}}{\sigma_{n+1}}}&\frac{i\kappa}{\sigma_{n+1}}e^{-i\frac{\kappa x_{n+1}}{\sigma_{n+1}}}&&\frac{-\sigma_{n}^2}{\sigma_{n+1}^2}e^{-i\frac{\kappa x_{n}}{\sigma_{n+1}}}&e^{-i\frac{\kappa x_{n+1}}{\sigma_{n+1}}}\\
\hdashline
\frac{i\kappa}{\sigma_1}e^{i\frac{\kappa x_0}{\sigma_1}}&\frac{-i\kappa}{\sigma_1}e^{i\frac{\kappa x_1}{\sigma_1}}&&\frac{-\sigma_0^2}{\sigma_1^2}e^{i\frac{\kappa x_0}{\sigma_1}}&e^{i\frac{\kappa x_1}{\sigma_1}}\\
\hspace{.7in}\ddots&\hspace{.7in}\ddots&&\hspace{.7in}\ddots&\hspace{.7in}\ddots\\
&\frac{i\kappa}{\sigma_{n+1}}e^{i\frac{\kappa x_{n+1}}{\sigma_{n+1}}}&\frac{-i\kappa}{\sigma_{n+1}}e^{i\frac{\kappa x_{n+1}}{\sigma_{n+1}}}&&\frac{-\sigma_{n}^2}{\sigma_{n+1}^2}e^{i\frac{\kappa x_{n}}{\sigma_{n+1}}}&e^{i\frac{\kappa x_{n+1}}{\sigma_{n+1}}}\\
&&\beta_3&&&\beta_4
\end{array}
\right).
\end{align*}}
\end{subequations}
The matrix $\mathcal{A}^F(\kappa)$ is made up of four $(n+2)\times(n+2)$ blocks as indicated by the dashed lines.  The two blocks in the upper half of $\mathcal{A}^F(\kappa)$ are zero except for entries on the main and $-1$ diagonals.  The lower two blocks of $\mathcal{A}^F(\kappa)$ only have entries on the main and $+1$ diagonals.  

As before we use Cramer's Rule to solve this system.  After multiplying the solutions by $\kappa e^{-\kappa^2t}$, integrating around $D_R^+$, and simplifying as in the previous section we follow a similar process to show the terms from $\mathcal{Y}^F(\kappa,T)$ do not contribute to our solution formula using Jordan's Lemma and Cauchy's Theorem.  One can show that $A_j^F(\kappa,T)$ can be replaced by $A_j^F(\kappa,t)$ by writing $\int_0^T \cdot \ud s$ as $\int_0^t \cdot \ud s+ \int_t^T \cdot \ud s$ and noticing where the function in analytic and decaying.  If the boundary conditions~\eqref{I2I_bcs} are time-independent then so is $A_j^F$.

In general, the initial-to-interface map for the heat equation on a finite domain with $n$ interfaces is given by
\begin{subequations}\label{I2I_INTg01solns_final_f}
\begin{align}
u^{(1)}(x_0,t) =&\int_{\partial D_R^+} e^{-\kappa^2t} \frac{\kappa\det(A^F_{1}(\kappa,t))}{i\pi\det(\mathcal{A}^F(\kappa))}\ud \kappa, \\
u^{(1)}_x(x_0,t)=&\int_{\partial D_R^+} e^{-\kappa^2t}\frac{\kappa \det(A^F_{n+3}(\kappa,t))}{i\pi\det(\mathcal{A}^F(\kappa))}\ud \kappa,\\
u^{(j)}(x_j,t) =&\int_{\partial D_R^+} e^{-\kappa^2t} \frac{\kappa\det(A^F_{j+1}(\kappa,t))}{i\pi\det(\mathcal{A}^F(\kappa))}\ud \kappa, \\
u^{(j)}_x(x_j,t)=&\int_{\partial D_R^+} e^{-\kappa^2t}\frac{\kappa \det(A^F_{j+n+3}(\kappa,t))}{i\pi\det(\mathcal{A}^F(\kappa))}\ud \kappa.
\end{align}
\end{subequations}
for $1\leq j\leq n+1$, where $A_j^F(\kappa,t)$ is the matrix $\mathcal{A}^F(\kappa,t)$ with the $j^\textrm{th}$ column replaced by $Y^F(\kappa,t)$.

As an example of a specific initial-to-interface map we consider~\eqref{n_heat} on a finite domain with $n=1$ with boundary conditions
\begin{align}
u^{(1)} _x(0 ,t)=&f_1(t) ,&t>0,\\
u^{(n+1)} _x(x_{2} ,t)=&f_2(t),&t>0.
\end{align}
and zero initial conditions
\begin{align}
u^{(j)} (x ,0)=&0 ,&x_{j-1}<x<x_j,
\end{align}
for $1\leq j\leq n+1$.  Using~\eqref{I2I_INTg01solns_final_f} we have
\begin{align*}
\sigma_1^2u^{(1)}_x(x_1,t)=&\frac{i\sigma_1\sigma_2}{\pi}\int_{\partial D_R^+}   \frac{\kappa e^{-\kappa^2t}\left(\sigma_1 \tilde{f_1}(\kappa^2,t)\sin\left(\frac{\kappa(x_1-x_2)}{\sigma_2}\right)-\sigma_2 \tilde{f_2}(\kappa^2,t)\sin\left(\frac{\kappa x_1}{\sigma_1}\right)\right)}{\sigma_1\cos\left(\frac{\kappa(x_1-x_2)}{\sigma_2}\right)\sin\left(\frac{\kappa x_1}{\sigma_1}\right)-\sigma_2\cos\left(\frac{\kappa x_1}{\sigma_1} \right)\sin\left(\frac{\kappa(x_1-x_2)}{\sigma_2}\right)} \ud \kappa,\\
u^{(1)}(x_1,t)=&\frac{-i}{\pi}\int_{\partial D_R^+}   \frac{e^{-\kappa^2t}\left(\sigma_1^2\tilde{f_1}(\kappa^2,t)\cos\left(\frac{\kappa(x_1-x_2)}{\sigma_2}\right)-\sigma_2^2\tilde{f_2}(\kappa^2,t)\cos\left(\frac{\kappa x_1}{\sigma_1}\right)\right)}{\sigma_1\cos\left(\frac{\kappa(x_1-x_2)}{\sigma_2}\right)\sin\left(\frac{\kappa x_1}{\sigma_1}\right)-\sigma_2\cos\left(\frac{\kappa x_1}{\sigma_1} \right)\sin\left(\frac{\kappa(x_1-x_2)}{\sigma_2}\right)} \ud \kappa.
\end{align*}

\section*{Acknowledgements}
N.E.S. acknowledges support from the National Science Foundation under grant number NSF-DGE-0718124.  Any opinions, findings, and conclusions or recommendations expressed in this material are those of the authors and do not necessarily reflect the views of the funding sources.

\bibliographystyle{plain}
\bibliography{FullBib}

\end{document}